\documentclass[reqno,dvips]{amsart}

\usepackage{mathrsfs}
\usepackage{amsmath}
\usepackage{amssymb}
\usepackage{cite}
\usepackage{latexsym}
\usepackage{graphicx}

\usepackage{color}

\theoremstyle{plain}
\newtheorem{thm}{Theorem}[section]

\newtheorem{dfn}[thm]{Definition}
\newtheorem{prop}[thm]{Proposition}
\newtheorem{rmk}[thm]{Remark}

\def\D{\mathrm{D}}

\def\I{\mathscr{I}}
\def\P{\mathscr{P}}
\def\R{\mathrm{R}}
\def\T{\mathscr{T}}

\def\d{\mathrm{d}}
\def\h{\mathrm{h}}
\def\t{\mathrm{T}}

\def\Cset{\mathbb{C}}

\def\Nset{\mathbb{N}}
\def\Qset{\mathbb{Q}}
\def\Rset{\mathbb{R}}
\def\Sset{\mathbb{S}}
\def\Tset{\mathbb{T}}
\def\Zset{\mathbb{Z}}

\def\epsilon{\varepsilon}

\DeclareMathOperator{\rank}{rank}
\DeclareMathOperator{\sn}{sn}
\DeclareMathOperator{\cn}{cn}
\DeclareMathOperator{\dn}{dn}

\DeclareMathOperator{\sech}{sech}
\DeclareMathOperator{\arccosh}{arccosh}

\makeatletter
 \@addtoreset{equation}{section}
\makeatother
\def\theequation{\arabic{section}.\arabic{equation}}

\begin{document}


\title[Nonintegrability of nearly integrable dynamical systems]%
{Nonintegrability of nearly integrable dynamical systems near resonant periodic orbits}

\author[K. Yagasaki]{Kazuyuki Yagasaki}

\address{Department of Applied Mathematics and Physics, Graduate School of Informatics,
Kyoto University, Yoshida-Honmachi, Sakyo-ku, Kyoto 606-8501, JAPAN}
\email{yagasaki@amp.i.kyoto-u.ac.jp}

\date{\today}

\begin{abstract}
In a recent paper by the author
 (K.~Yagasaki, Nonintegrability of the restricted three-body problem, submitted for publication),
 a technique was developed for determining
 whether nearly integrable systems are not meromorphically Bogoyavlenskij-integrable
 such that the first integrals and commutative vector fields also depend meromorphically
 on the small parameter.
Here we continue to demonstrate the technique for some classes of dynamical systems.
In particular, we consider time-periodic perturbations of single-degree-of-freedom Hamiltonian systems
 and discuss a relationship of the technique with the subharmonic Melnikov method,
 which enables us to detect the existence of periodic orbits and their stability.
We illustrate the theory for the periodically forced Duffing oscillator
 and two more additional examples: second-order coupled oscillators
 and  a two-dimensional system of pendulum-type subjected to a constant torque.
\end{abstract}
\maketitle


\section{Introduction}

In this paper we consider systems of the form 
\begin{equation}
\dot{I}=\epsilon h(I,\theta;\epsilon),\quad
\dot{\theta}=\omega(I)+\epsilon g(I,\theta;\epsilon),\quad
(I,\theta)\in\Rset^\ell\times\Tset^m,
\label{eqn:aasys}
\end{equation}
and study its nonintegrability near resonant periodic orbits,
 where $\ell,m\in\Nset$,  $\Tset^m=(\Rset/2\pi\Zset)^m$,
 $\epsilon$ is a small parameter such that $0<|\epsilon|\ll 1$,
 and $\omega:\Rset^\ell\to\Rset^m$, $h:\Rset^\ell\times\Tset^m\times\Rset\to\Rset^\ell$
 and $g:\Rset^\ell\times\Tset^m\times\Rset\to\Rset^m$
 are meromorphic or analytic in the arguments.
We extend the domain of the independent variable $t$
 to a domain including $\Rset$ in $\Cset$ and do so for the dependent variables.
The system~\eqref{eqn:aasys} is Hamiltonian if $\ell=m$ as well as $\epsilon=0$ or
\[
\D_I h(I,\theta;\epsilon)\equiv-\D_\theta g(I,\theta;\epsilon),
\]
and non-Hamiltonian if not.
When $\epsilon=0$, Eq.~\eqref{eqn:aasys} becomes
\begin{equation}
\dot{I}=0,\quad
\dot{\theta}=\omega(I)
\label{eqn:aasys0}
\end{equation}
which we refer to as the \emph{unperturbed system} for \eqref{eqn:aasys}.
Here we adopt the following definition of integrability due to Bogoyavlenskij  \cite{B98}.

\begin{dfn}[Bogoyavlenskij]
\label{dfn:1a}
For $n\in\Nset$ an $n$-dimensional dynamical system
\[
\dot{x}=f(x),\quad x\in\Rset^n\text{ or }\Cset^n,
\]
is called \emph{$(q,n-q)$-integrable} or simply \emph{integrable} 
 if there exist $q$ vector fields $f_1(x)(:= f(x)),f_2(x),\dots,f_q(x)$
 and $n-q$ scalar-valued functions $F_1(x),\dots,F_{n-q}(x)$ such that
 the following two conditions hold:
\begin{enumerate}
\setlength{\leftskip}{-1.8em}
\item[\rm(i)]
$f_1(x),\dots,f_q(x)$ are linearly independent almost everywhere and commute with each other,
 i.e., $[f_j,f_k](x):=\D f_k(x)f_j(x)-\D f_j(x)f_k(x)\equiv 0$ for $j,k=1,\ldots,q$,
 where $[\cdot,\cdot]$ denotes the Lie bracket$;$
\item[\rm(ii)]
The derivatives $\D F_1(x),\dots, \D F_{n-q}(x)$ are linearly independent almost everywhere
 and $F_1(x),\dots,F_{n-q}(x)$ are first integrals of $f_1, \dots,f_q$,
 i.e., $\D F_k(x)\cdot f_j(x)\equiv 0$ for $j=1,\ldots,q$ and $k=1,\ldots,n-q$,
 where ``$\cdot$'' represents the inner product.
\end{enumerate}
We say that the system is \emph{meromorphically} $($resp. \emph{analytically}$)$ \emph{integrable}
 if the first integrals and commutative vector fields are meromorphic $($resp. analytic$)$. 
\end{dfn}

Definition~\ref{dfn:1a} is considered as a generalization of 
 Liouville-integrability for Hamiltonian systems \cite{A89,M99}
 since an $n$-degree-of-freedom Liouville-integrable Hamiltonian system with $n\ge 1$
 has not only $n$ functionally independent first integrals
 but also $n$ linearly independent commutative (Hamiltonian) vector fields
 generated by the first integrals.
The unperturbed system~\eqref{eqn:aasys0} is meromorphically or analytically $(m,\ell)$-integrable
 in the Bogoyavlenskij sense:
 $F_j(I,\theta)=I_j$, $j=1,\ldots,\ell$, are first integrals
 and $f_j(I,\theta)=(0,e_j)\in\Rset^\ell\times\Rset^m$, $j=2,\ldots,m$,
 give $m$ commutative vector fields along with its own vector field,
 where $e_j$ is the $m$-dimensional vector of which the $j$th element is the unit
 and the other elements are zero. 
Conversely, a general $(m,\ell)$-integrable system is transformed to the form \eqref{eqn:aasys0}
 if the level set for the first integrals $F_1(x),\ldots,F_m(x)$ has a connected compact component.
See \cite{B98,MY21b,Z18} for more details.
Thus, the system \eqref{eqn:aasys} can be regarded
 as a normal form for perturbations of general $(m,\ell)$-integrable systems.

In a recent paper \cite{Y21a},
 a technique was developed for determining
 whether the system~\eqref{eqn:aasys} is not meromorphically Bogoyavlenskij-integrable
 such that the first integrals and commutative vector fields also depend meromorphically
 on the small parameter $\epsilon$ near $\epsilon=0$.
Moreover, the technique was applied to prove that
 the restricted three-body problem are not meromorphically integrable
 in both the planar and spatial cases
 even if the first integrals are not required to depend meromorphically on the parameter,
 the mass ratio of the primaries.
The basic idea used there was similar to that of Morales-Ruiz \cite{M02},
 who studied time-periodic Hamiltonian perturbations
 of single-degree-of-freedom Hamiltonian systems
 and showed a relationship of their nonintegrability
 with a version due to Ziglin \cite{Z82} of the Melnikov method \cite{M63}.
The Melnikov method enables us
 to detect transversal self-intersection of complex separatrices of periodic orbits
 unlike the standard version \cite{GH83,M63,W90}.
More concretely, under some restrictive conditions,
 he essentially proved that they are meromorphically nonintegrable
 when the small parameter is taken as one of the state variables
 if the Melnikov functions are not identically zero,
 based on a generalized version due to Ayoul and Zung \cite{AZ10}
 of the Morales-Ramis theory \cite{M99,MR01}.
Their generalized versions for the Morales-Ramis theory and its extension,
 the Morales-Ramis-Sim\'o theory \cite{MRS07}, were also used in \cite{Y21a}.
The developed technique was also applied
 to give a new proof of Poincar\'e's result of \cite{P92} on the restricted three-body problem
 in \cite{Y21b}.

In this paper we continue to demonstrate the technique of \cite{Y21a}
 for some classes of dynamical systems.
In particular, we consider time-periodic perturbations of single-degree-of-freedom Hamiltonian systems
 and discuss a relationship of the technique
 with the subharmonic Melnikov method \cite{GH83,W90,Y96},
 which enables us to detect the existence of periodic orbits and their stability and bifurcations,
 like Morales-Ruiz \cite{M02} for homoclinic orbits.
So we show that they are nonintegrable in the meaning stated above
 if certain complex integrals similar to the subharmonic Melnikov functions are not zero.
See Theorem~\ref{thm:3a} below for the precise statement.
The similarity of this result to that of \cite{M02} is very remarkable.

We also illustrate the theory for the periodically forced Duffing oscillator
\[
\ddot{w}+\epsilon\delta\dot{w}+aw+w^3=\epsilon\beta\cos\nu t,\quad
w\in\Rset,
\]
or as a first-order system
\begin{equation}
\dot{x}_1=x_2,\quad
\dot{x}_2=-x_1-x_1^3+\epsilon(\beta\cos\nu t-\delta x_2),\quad
x_1,x_2\in\Rset,
\label{eqn:duf}
\end{equation}
where $a=\pm1$ or $0$, and $\beta,\nu>0$ and $\delta\ge 0$ are constants.
It is well-known that Duffing \cite{D18} studied this type of system early in the twentieth century 
 but it is interesting that Poincar\'e also discussed the existence of periodic solutions
 for $\delta=0$ about the end of the nineteenth century in his memoir \cite{P90}.
See also Section~5.6 of \cite{B96}.
Holmes \cite{H79} used the homoclinic Melnikov method \cite{GH83,M63,W90}
 to prove the occurrence of transverse intersection
 between the stable and unstable manifolds of a periodic orbit near $(x_1,x_2)=(0,0)$ 
 for $a=-1$ with $\epsilon>0$ sufficiently small.
The occurrence of such transverse intersection implies, e.g., by Theorem~3.10 of \cite{M73},
 the real-analytic nonintegrability near the unperturbed homoclinic orbit.
Motonaga and Yagasaki \cite{MY21b}
 showed the real-analytic nonintegrability of \eqref{eqn:duf} with $a=-1$
 near the unperturbed homoclinic orbits in the meaning stated above
 even when such transverse intersection does not occur
 (see Remark~\ref{rmk:4c}(ii) for more details).
Ueda  \cite{U78} also found chaotic motions in both analog and numerical simulations
 when $a=0$ but $\epsilon$ is not small.
Moreover, the rational nonintegrability of the parametric excitation case, e.g.,
\[
\dot{x}_1=x_2,\quad
\dot{x}_2=ax_1-x_1^3-\delta x_2+\beta x_1\cos\nu t,
\]
was recently proved in \cite{MY18} when $\mathrm{e}^{i\nu t}=\cos\nu t+i\sin\nu t$
 is taken as a state variable.
So the Duffing oscillator \eqref{eqn:duf} has been believed to be nonintegrable
 besides near the unperturbed homoclinic orbits for $a=-1$,
 but its proof has not been given.
We show that the system \eqref{eqn:duf} is meromorphically nonintegrable
 near the resonant periodic orbits in the meaning stated above
 when $a=\pm 1$ and $0$.

Moreover, we give two more concrete examples.
The first one is second-order coupled oscillators
 of which the special case is often referred
 to as the \emph{second-order Kuramoto model} \cite{RPJK16}.
The second one is a two-dimensional system of pendulum-type subjected to a constant torque.
We show that it is not integrable as a system on $\Cset\times(\Cset/2\pi\Zset)$
 although it has a first integral as a system on $\Rset^2$ or $\Cset^2$.

This paper is organized as follows:
In Section~2 we review the technique of \cite{Y21a} in a necessary context.
In Section~3 we apply the technique
 to time-periodic perturbations of single-degree-of-freedom Hamiltonian systems
 and discuss a relationship of the result with the subharmonic Melnikov method.
We illustrate the theory for the periodically forced Duffing oscillator \eqref{eqn:duf} in Section~4.
Finally, we provide the additional two examples in Section~5.


\section{General Technique}

In this section
 we review the technique of \cite{Y21a} for the nonintegrability of \eqref{eqn:aasys}.
We make the following assumption on the unperturbed system \eqref{eqn:aasys0}:
\begin{enumerate}
\setlength{\leftskip}{-1em}
\item[\bf(A1)]
For some $I^\ast\in\Rset^\ell$, a resonance of multiplicity $m-1$,
\[
\dim_\Qset\langle\omega_1(I^\ast),\ldots,\omega_m(I^\ast)\rangle=1,
\]
occurs with $\omega(I^\ast)\neq 0$,
 i.e., there exists a constant $\omega^\ast>0$ such that 
\[
\frac{\omega(I^\ast)}{\omega^\ast}\in\Zset^m\setminus\{0\},
\]
where $\omega_j(I)$ is the $j$th element of $\omega(I)$ for $j=1,\ldots,m$.
\end{enumerate}
Note that we can replace $\omega^\ast$ with $\omega^\ast/k$ for any $k\in\Nset$ in (A1).
We refer to the $m$-dimensional torus $\T^\ast=\{(I^\ast,\theta)\mid\theta\in\Tset^m\}$
 as the \emph{resonant torus}
 and to periodic orbits $(I,\theta)=(I^\ast,\omega(I^\ast)t+\theta_0)$, $\theta_0\in\Tset^m$,
 on $\T^\ast$ as the \emph{resonant periodic orbits}.
Let $T^\ast=2\pi/\omega^\ast$.
We also make the following assumption.

\begin{figure}
\includegraphics[scale=0.8
]{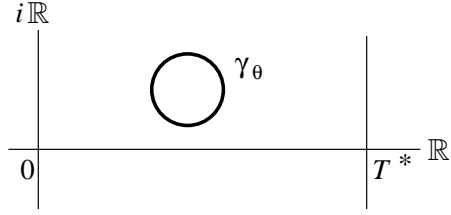}
\caption{Assumption~(A2).\label{fig:1a}}
\end{figure}

\begin{enumerate}
\setlength{\leftskip}{-1em}
\item[\bf(A2)]
For some $\theta\in\Tset^m$,
 there exists a closed loop $\gamma_\theta$
 in a domain including $(0,T^\ast)\subset\Rset$ in $\Cset$
 such that $\gamma_\theta\cap(i\Rset\cup(T^\ast+i\Rset))=\emptyset$ and
\begin{equation}
\I(\theta)
:=\D\omega(I^\ast)\int_{\gamma_\theta}h(I^\ast,\omega(I^\ast)\tau+\theta;0)\d\tau
\label{eqn:A2}
\end{equation}
is not zero. 
See Fig.~\ref{fig:1a}
\end{enumerate}
Note that the condition $\gamma_\theta\cap(i\Rset\cup(T^\ast+i\Rset))=\emptyset$
 is not essential in (A2), since it always holds
 by replacing $\omega^\ast$ with $\omega^\ast/k$ for sufficiently large $k\in\Nset$ if necessary.
We can prove the following theorem which guarantees that
 conditions~(A1) and (A2) are sufficient for nonintegrability of \eqref{eqn:aasys}
 in the meaning stated in Section~1.

\begin{thm}
\label{thm:main}
Let $\Gamma$ be any domain in $\Cset/T^\ast\Zset$ containing $\Rset/T^\ast\Zset$ and $\gamma_\theta$.
Suppose that assumption~{\rm(A1)} and {\rm(A2)} hold for some $\theta_0\in\Tset^m$.
Then the system \eqref{eqn:aasys} is not meromorphically integrable in the Bogoyavlenskij sense
 near the resonant periodic orbit $(I,\theta)=(I^\ast,\omega(I^\ast)\tau+\theta_0)$
 with $\tau\in\Gamma$ 
 such that the first integrals and commutative vector fields also depend meromorphically
 on $\epsilon$ near $\epsilon=0$.
Moreover, if {\rm(A2)} holds for $\theta\in\Delta$, where $\Delta$ is a dense set in $\Tset^m$,
 then the conclusion holds for any resonant periodic orbit on the resonant torus $\T^\ast$.
\end{thm}

See Section~2 of \cite{Y21a} for a proof of Theorem~\ref{thm:main}.
A more general result was obtained there.

Systems of the form \eqref{eqn:aasys} have attracted much attention,
 especially when they are Hamiltonian.
See \cite{A89,AKN06,K96} and references therein for more details.
In particular, Kozlov \cite{K96} extended the famous result of Poincar\'e \cite{P90,P92}
 for Hamiltonian systems to the general case of \eqref{eqn:aasys}
 and gave sufficient conditions for nonexistence of additional real-analytic first integrals
 depending analytically on $\epsilon$ near $\epsilon=0$.
See also \cite{AKN06,K83} for his result in Hamiltonian systems. 
Moreover, Motonaga and Yagasaki \cite{MY21b} gave sufficient conditions
 for the system \eqref{eqn:aasys} to be real-analytically nonintegrable in the Bogoyavlenskij sense
 such that the first integrals and commutative vector fields also depend real-analytically
 on $\epsilon$ near $\epsilon=0$.
Some details on these results are provided in our context
 and compared with Theorem~\ref{thm:main} in Appendix~A.
We remark that
 the results of \cite{K96,P92,MY21b} say nothing about the integrability of \eqref{eqn:aasys}
 under the hypotheses of Theorem~\ref{thm:main}.


\section{Time-Periodic Perturbations of Single-Degree-of-Freedom Hamiltonian Systems}

We next apply the technique of Section~2
 to time-periodic perturbations of single-degree-of-freedom Hamiltonian systems,
 and discuss a relationship of our result with the subharmonic Melnikov method \cite{GH83,W90,Y96},
 as in the related work \cite{MY21a,MY21b}.
 
\begin{figure}
\includegraphics[scale=1
]{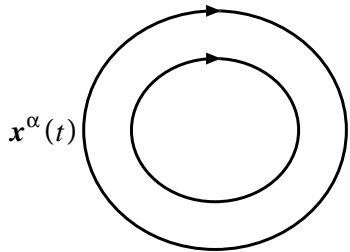}
\caption{Assumption~(M1).\label{fig:3a}}
\end{figure}

Consider two-dimensional systems of the form 
\begin{equation}
\dot{x}=J\D H(x)+\epsilon u(x,\nu t),\quad
x\in\Rset^2,
\label{eqn:syse}
\end{equation}
where $\nu>0$ is a constant, $H:\Rset^2\to\Rset$ and $u:\Rset^2\times\Sset^1$ are analytic,
 and $J$ is the $2\times 2$ symplectic matrix,
\[
J=\begin{pmatrix}
0 & 1\\
-1 & 0
\end{pmatrix}.
\]
When $\epsilon=0$, Eq.~\eqref{eqn:syse} becomes a planar Hamiltonian system
\begin{equation}
\dot{x}=J\D_x H(x)
\label{eqn:sys0}
\end{equation}
with a Hamiltonian function $H(x)$.
We make the following assumptions on the unperturbed system \eqref{eqn:sys0}:
\begin{enumerate}
\setlength{\leftskip}{-0.8em}
\item[\bf(M1)]
There exists a one-parameter family of periodic orbits $x^\alpha(t)$, $\alpha\in(\alpha_1,\alpha_2)$,
 with period $T^\alpha>0$ for some $\alpha_1<\alpha_2$ (see Fig.~\ref{fig:3a});
\item[\bf(M2)]
$x^\alpha(t)$ is analytic with respect to $\alpha\in(\alpha_1,\alpha_2)$.
\end{enumerate}
Note that in assumption~(M1)
 $x^\alpha(t)$ is automatically analytic with respect to $t$
 since the vector field of \eqref{eqn:sys0} is analytic.
Following an approach of \cite{Y96},
 we can transform \eqref{eqn:syse} into the form \eqref{eqn:aasys} as follows.

We first define the scalar action variable $I^\alpha$
 for each periodic orbit $x^\alpha(t)=(x_1^\alpha(t),x_2^\alpha(t))$ as
\begin{equation}
I^\alpha=\frac{1}{2\pi}\int_{x^\alpha}x_2\d x_1
 =\frac{1}{2\pi}\int_0^{T^\alpha}x_2^\alpha(t)\dot{x}_1^\alpha(t)\d t
\label{eqn:I1}
\end{equation}
in the standard manner (see, e.g., Chapter~10 of \cite{A89}).
The action variable $I$ can thus be determined only by $\alpha$.
We assume that $\d\alpha/\d I>0$ without loss of generality,
 and apply the implicit function theorem to \eqref{eqn:I1} to represent $\alpha$ as a function of $I$:
 $\alpha=\alpha(I)$.
We can show that the symplectic transformation from $(I,\theta_1)$ to $x$ is given by
\begin{equation}
x=x^{\alpha(I)}\left(\frac{\theta_1}{\Omega(I)}\right),
\label{eqn:trans}
\end{equation}
where
\[
\Omega(I)=\frac{2\pi}{T^{\alpha(I)}}.
\]
We see that $\d\Omega/\d I\neq 0$ at $I=I^\alpha$ if $\d T^\alpha/\d\alpha\neq 0$.
Moreover, we have the relations
\begin{equation}
\D_xI=-J\frac{\partial x}{\partial\theta_1},\quad
\D_x\theta_1=J\frac{\partial x}{\partial I}.
\label{eqn:rel}
\end{equation}
Let $\theta_2=\nu t$ in \eqref{eqn:syse}.
Using \eqref{eqn:syse}, \eqref{eqn:trans} and \eqref{eqn:rel}, we obtain
\begin{equation}
\dot{I}=\epsilon h(I,\theta_1,\theta_2),\quad
\dot{\theta}_1=\Omega(I)+\epsilon g_1(I,\theta_1,\theta_2),\quad
\dot{\theta}_2=\nu,
\label{eqn:aasys4}
\end{equation}
where
\begin{align*}
&
h(I,\theta_1,\theta_2)
 =\frac{1}{\Omega(I)}\D H\left(x^{\alpha(I)}\left(\frac{\theta_1}{\Omega(I)}\right)\right)
 \cdot u\left(x^{\alpha(I)}\left(\frac{\theta_1}{\Omega(I)}\right),\theta_2\right),\\
&
g_1(I,\theta_1,\theta_2)
 =J\frac{\partial}{\partial I}x^{\alpha(I)}\left(\frac{\theta_1}{\Omega(I)}\right)
 \cdot u\left(x^{\alpha(I)}\left(\frac{\theta_1}{\Omega(I)}\right),\theta_2\right).
\end{align*}
See Section~2 of \cite{Y96} for the details on these computations.
The system \eqref{eqn:aasys4} has the form \eqref{eqn:aasys}
 with $\ell=1$, $m=2$ and $\omega(I)=(\Omega(I),\nu)^\t$,
  where the superscript ``$\t$'' represents the transpose operator.

We assume that at $\alpha=\alpha^{l/n}$
\begin{equation*}
\frac{2\pi}{T^\alpha}=\frac{n}{l}\nu,
\end{equation*}
where $l$ and $n$ are relatively prime integers,
 so that assumption~(A1) holds with $\omega^\ast=2\pi/nT^\alpha=\nu/l$.
We define the \emph{subharmonic Melnikov function} as
\begin{equation}
M^{l/n}(\phi)=\int_0^{2\pi l/\nu}\D H(x^{\alpha}(t))\cdot u(x^\alpha(t),\nu t+\phi)\d t,
\label{eqn:mel}
\end{equation}
where $\alpha=\alpha^{l/n}$.
If $M^{l/n}(\phi)$ has a simple zero at $\phi=\phi_0$ and $\d T^\alpha/\d\alpha\neq 0$,
 i.e., $\d\Omega(I^\alpha)/\d I\neq 0$,
 then there exists a periodic orbit near $(x,\phi)=(x^\alpha(t),\nu t+\phi_0)$ in \eqref{eqn:syse}.
See Theorem~3.1 of \cite{Y96}.
A similar result is also found in \cite{GH83,W90}.
The stability of the periodic orbit can also be determined easily \cite{Y96}.
Moreover, several bifurcations of periodic orbits
 when $\d\Omega(I^\alpha)/\d I\neq 0$ or not
 were discussed in \cite{Y96,Y02,Y03}.

Noting that $\Omega(I^\alpha)=n\nu/l$ at $\alpha=\alpha^{l/n}$
 and applying Theorem~\ref{thm:main} to \eqref{eqn:aasys4},
 we obtain the following.
 
\begin{thm}
\label{thm:3a}
Suppose that at $\alpha=\alpha^{l/n}$, $\d T^\alpha/\d\alpha\neq 0$
 and there exists a closed loop $\gamma_\phi$
 in a domain including $(0,2\pi l/\nu)$ in $\Cset$
 such that $\gamma_\phi\cap(i\Rset\cup(2\pi l/\nu+i\Rset))=\emptyset$ and 
\begin{equation}
\hat{\I}(\phi)
=\int_{\gamma_\phi}\D H(x^\alpha(\tau))
 \cdot u\left(x^\alpha(\tau),\nu\tau+\phi\right)\d\tau
\label{eqn:thm4a}
\end{equation}
is not zero for some $\phi=\phi_0\in\Sset^1$.
Then the system \eqref{eqn:aasys4}, equivalently \eqref{eqn:syse},
 is not meromorphically integrable in the meaning of Theorem~$\ref{thm:main}$
 near the resonant periodic orbit $(x,\phi)=(x^\alpha(t),\nu t+\phi_0)$ with $\alpha=\alpha^{l/n}$
 on any domain $\hat{\Gamma}$ in $\Cset/(2\pi l/\nu)\Zset$ containing $\Rset/(2\pi l/\nu)\Zset$ and $\gamma_\phi$.
Moreover, if the integral $\hat{\I}(\phi)$ is not zero for any $\phi\in\hat{\Delta}$,
 where $\hat{\Delta}$ is a dense set of $\Sset^1$,
 then the conclusion holds for any periodic orbit on the resonant torus
 $\T^\ast=\{(x^\alpha(\tau),\nu\tau+\phi)\mid\tau\in\hat{\Gamma},\phi\in\Sset^1,\alpha=\alpha^{l/n}\}$.
\end{thm}

\begin{rmk}
\label{rmk:3a}
Let $U$ be a neighborhood of $\alpha=\alpha_0\in (\alpha_1,\alpha_2)$.
From Theorem~$\ref{thm:MY}$ we obtain the following for \eqref{eqn:syse}
 $($see Theorem~$5.2$ of {\rm\cite{MY21b}):}
If there exists a key set $D\subset D_\R:=\{\alpha^{l/n}\in U\mid\text{$l,n\in\Nset$ are relatively prime}\}$ for $C^\omega(U)$
 such that $M^{l/n}(\phi)$ is not constant for $\alpha^{l/n}\in D$,
 then for $|\epsilon|\neq 0$ sufficiently small
 the system \eqref{eqn:syse} is not real-analytically integrable
 in the meaning of Theorem~$\ref{thm:MY}$
 near $\{x^\alpha(t)\mid t\in[0,T^\alpha)\}\times\Sset^1$
 with $\alpha=\alpha_0$.
Note that $D_\R$ is a key set for $C^\omega(U)$.
\end{rmk}

Note that the integrand in \eqref{eqn:thm4a} is the same
 as in the Melnikov function \eqref{eqn:mel} although the path of integration is different.
An integral similar to \eqref{eqn:thm4a} for not periodic but homoclinic orbits was used
 in \cite{M02,Z82}.


\section{Periodically forced Duffing oscillator}

\begin{figure}[t]
\includegraphics[scale=0.28
]{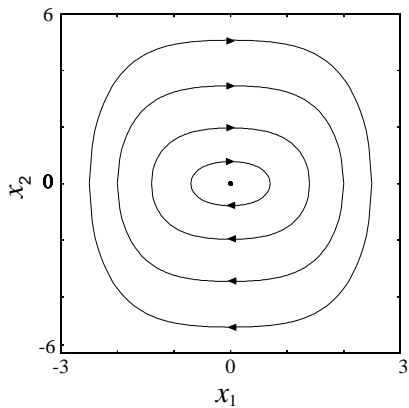}
\caption{Phase portraits of \eqref{eqn:duf} with $\epsilon=0$ and $a=1$.
\label{fig:4a}}
\end{figure}

We now consider the periodically forced Duffing oscillator \eqref{eqn:duf}
 and apply Theorem~\ref{thm:3a}.
When $\epsilon=0$,
 Eq.~\eqref{eqn:duf} becomes a single-degree-of-freedom Hamiltonian system
 with the Hamiltonian
\[
H=\frac{1}{2}ax_1^2+\frac{1}{4}x_1^4+\frac{1}{2}x_2^2,
\]
and it is a special case of \eqref{eqn:sys0}.

\subsection{Case of $a=1$}
We begin with the case of $a=1$.
The phase portraits of \eqref{eqn:duf} with $\epsilon=0$ are shown in Fig.~\ref{fig:4a}.
In particular,
 there exists a one-parameter family of periodic orbits
\begin{align*}
x^k(t)
 =&\biggl(\frac{\sqrt{2}k}{\sqrt{1-2k^2}}\cn\left(\frac{t}{\sqrt{1-2k^2}}\right),\\
& \quad-\frac{\sqrt{2}k}{1-2k^2}\sn\left(\frac{t}{\sqrt{1-2k^2}}\right)
\dn\left(\frac{t}{\sqrt{1-2k^2}}\right)\biggr),\quad
k\in\bigl(0,1/\sqrt{2}\bigr),
\end{align*}
and their period is given by $T^k=4K(k)\sqrt{1-2k^2}$ (see \cite{Y94,Y96}),
 where sn, cn and dn represent the Jacobi elliptic functions, $k$ is the elliptic modulus
 and $K(k)$ is the complete elliptic integral of the first kind.
See, e.g., \cite{BF54,WW27} for general information on elliptic functions.
Assume that the resonance condition
\begin{equation}
nT^k=\frac{2\pi l}{\nu},\quad\mbox{i.e.,}\quad
\nu=\frac{\pi l}{2nK(k)\sqrt{1-2k^2}},
\label{eqn:res4a}
\end{equation}
holds at $k=k^{l/n}$ for $l,n>0$ relatively prime integers.
We compute the subharmonic Melnikov function \eqref{eqn:mel} for $x^k(t)$ as
\begin{align*}
M^{l/n}(\phi)=&\int_0^{2\pi l/\nu}x_2^k(t)(-\delta x_2^k(t)+\beta\cos(\nu t+\phi))\d t\notag\\
=&-\delta J_1(k,n)+\beta J_2(k,l,n)\sin\phi,
\end{align*}
where
\begin{align*}
&
J_1(k,n)=\frac{8n[(2k^2-1)E(k)+k'^2K(k)]}{3(1-2k^2)^{3/2}},\\
&
J_2(k,l,n)=
\begin{cases}
\displaystyle
2\sqrt{2}\pi\nu
 \sech\left(\frac{\pi lK(k')}{2K(k)}\right) & \mbox{(for $n=1$ and $l$ odd)};\\
0\quad & \mbox{(for $n\neq 1$ or $l$ even)}.
\end{cases}
\end{align*}
Here $E(k)$ is the complete elliptic integral of the second kind
 and $k'=\sqrt{1-k^2}$ is the complementary elliptic modulus.
See also \cite{Y94,Y96} for the computations of the Melnikov function.

On the other hand, we write the integral \eqref{eqn:thm4a} as
\begin{align}
\hat{\I}(\phi)=&
-\frac{2k^2 \delta}{(1-2k^2)^2}\int_{\gamma_\phi}
 \sn^2\left(\frac{\tau}{\sqrt{1-2k^2}}\right)
 \dn^2\left(\frac{\tau}{\sqrt{1-2k^2}}\right)\d\tau\notag\\
&
-\frac{\sqrt{2}k\beta}{1-2k^2}\int_{\gamma_\phi}
 \sn\left(\frac{\tau}{\sqrt{1-2k^2}}\right)
  \dn\left(\frac{\tau}{\sqrt{1-2k^2}}\right)\cos(\nu\tau+\phi)\d\tau.
\label{eqn:b1}
\end{align}
Letting $\gamma_\phi$ be a circle centered at $\tau=i\sqrt{1-2k^2}K(k')$ with sufficiently small radius,
 we compute
\begin{align}
\hat{\I}(\phi)
=-2\sqrt{2}\pi\nu\beta\biggl(\cosh\biggl(\frac{\pi l K(k')}{2nK(k)}\biggr)\sin\phi
  &-i\sinh\biggl(\frac{\pi l K(k')}{2nK(k)}\biggr)\cos\phi\biggr),
 \label{eqn:b2}
\end{align}
which is not zero for any $\phi\in\Sset^1$.
See Appendix~B  for the derivation of \eqref{eqn:b2}.
Applying Theorem~\ref{thm:3a}, we obtain the following.

\begin{prop}
\label{prop:4a}
Let $\hat{\Gamma}$ be a domain in $\Cset/(2\pi l/\nu)\Zset$
 containing $\Rset/(2\pi l/\nu)\Zset$ and $\tau=i\sqrt{1-2k^2}K(k')$.
The periodically forced Duffing oscillator \eqref{eqn:duf} with $a=1$ 
 is meromorphically nonintegrable in the meaning of Theorem~$\ref{thm:main}$
 near any periodic orbit on the resonant torus
 $\T^k=\{(x^k(\tau),\nu\tau+\theta)\mid\tau\in\hat{\Gamma},\theta\in\Sset^1,k=k^{l/n}\}$
 for $l,n>0$ relatively prime integers.
\end{prop}

\begin{rmk}\
\label{rmk:4a}
\begin{itemize}
\setlength{\leftskip}{-1.8em}
\item[(i)]
If $\beta=0$, then Proposition~$\ref{prop:4a}$ says nothing
 about the nonintegrability of \eqref{eqn:duf}
 since the integral \eqref{eqn:b2} is identically zero.

\item[(ii)]
For any neighborhood $U$ of $k\in(0,1/\sqrt{2})$
 there is not a key set $D\subset U$ for $C^\omega(U)$ such that $M^{l/n}(\phi)$ is not constant
 for $k\in D$ satisfying \eqref{eqn:res4a}.
Hence, Theorem~$\ref{thm:MY}$ is not applicable.
See Remark~$\ref{rmk:3a}$.
\end{itemize}
\end{rmk}

\subsection{Case of $a=0$}

\begin{figure}[t]
\includegraphics[scale=0.28
]{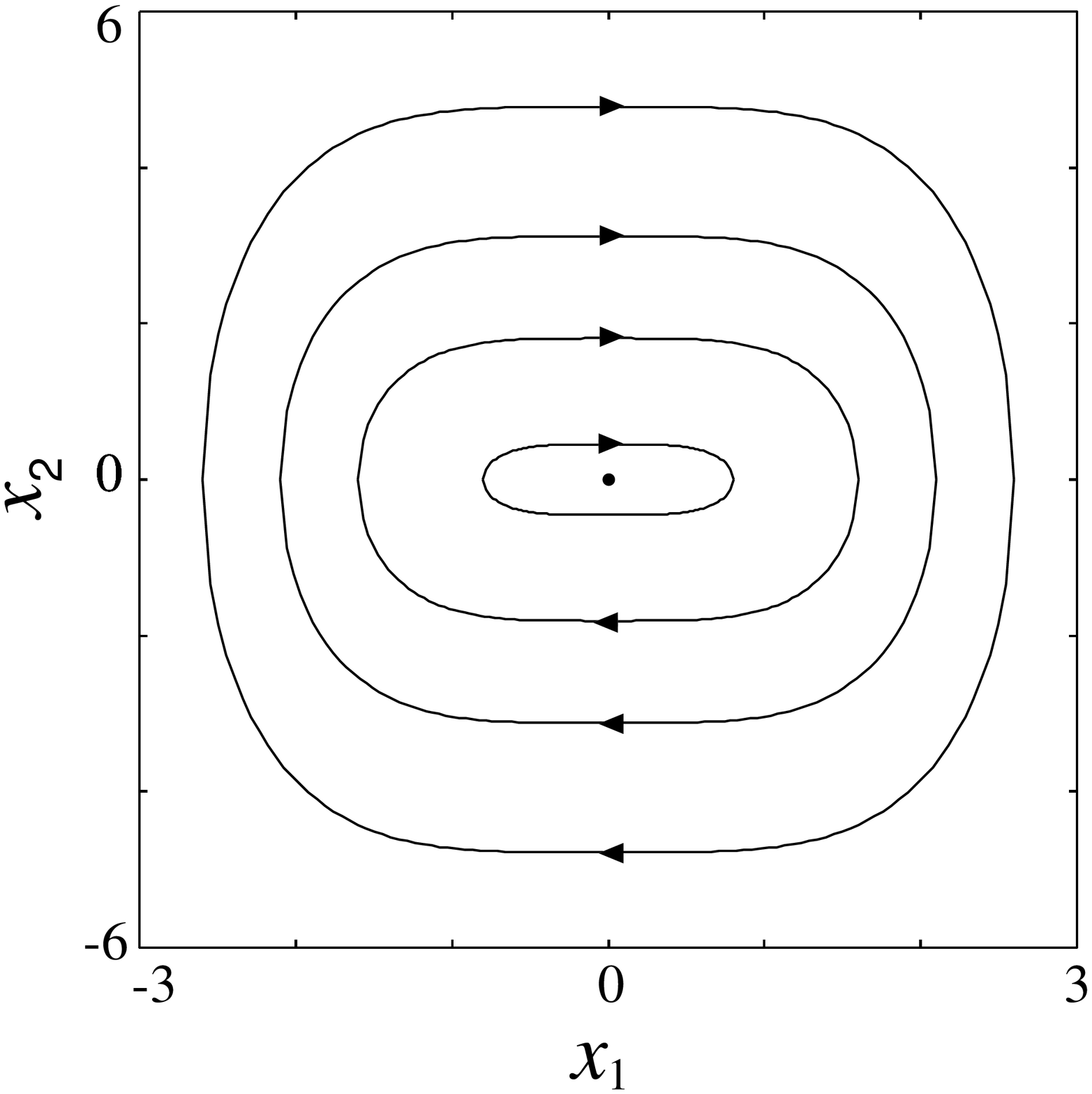}
\caption{Phase portraits of \eqref{eqn:duf} with $\epsilon=0$ and $a=0$.
\label{fig:4b}}
\end{figure}

We turn to the case of $a=0$ in \eqref{eqn:duf}.
The phase portraits of \eqref{eqn:duf} with $\epsilon=0$ are shown in Fig.~\ref{fig:4b}.
In particular, there exists a one-parameter family of periodic orbits
\begin{align*}
x^\alpha(t)
 =&(\alpha\cn \alpha t, -\alpha^2\sn\alpha t \dn \alpha t),\quad
\alpha\in\bigl(0,\infty),
\end{align*}
and their period is given by $T^\alpha=4K(1/\sqrt{2})/\alpha$,
 where the elliptic modulus in the Jacobi elliptic functions is $k=1/\sqrt{2}$
 and $K(1/\sqrt{2})=1.854\ldots$.
Assume that the resonance condition
\begin{equation}
nT^\alpha=\frac{2\pi l}{\nu},\quad\mbox{i.e.,}\quad
\nu=\frac{\pi l\alpha}{2nK(1/\sqrt{2})},
\label{eqn:res4b}
\end{equation}
holds at $\alpha=\alpha^{l/n}$ for $l,n>0$ relatively prime integers.
As in the case of $a=1$,
 we compute the subharmonic Melnikov function \eqref{eqn:mel} for $x^k(t)$ as
\begin{align*}
M^{l/n}(\phi)=&\int_0^{2\pi l/\nu}x_2^\alpha(t)(-\delta x_2^\alpha(t)+\beta\cos(\nu t+\phi))\d t\notag\\
=&-\delta J_1(\alpha,n)+\beta J_2(\alpha,l,n)\sin\phi,
\end{align*}
where
\begin{align*}
&
J_1(\alpha,n)=\frac{4n\alpha^3 K(1/\sqrt{2})}{3},\\
&
J_2(\alpha,l,n)=
\begin{cases}
\displaystyle
2\sqrt{2}\pi\nu
\sech\left(\frac{\pi l}{2}\right) & \mbox{(for $n=1$ and $l$ odd)};\\
0\quad & \mbox{(for $n\neq 1$ or $l$ even)}.
\end{cases}
\end{align*}
On the other hand, we write the integral \eqref{eqn:thm4a} as
\[
\hat{\I}(\phi)=-\alpha^4\delta\int_{\gamma_\phi}
 \sn^2\alpha t\,\dn^2\alpha\tau\d\tau
-\alpha^2\beta\int_{\gamma_\phi}
 \sn\alpha\tau\,\dn\alpha\tau\,\cos(\nu\tau+\phi)\d\tau.
\]
We take a circle centered at $\tau=i\alpha K(1/\sqrt{2})$ with sufficiently small radius as $\gamma_\phi$,
 and compute
\begin{equation}
\hat{\I}(\phi)
=-2\sqrt{2}\pi\nu\beta
 \biggl(\cosh\biggl(\frac{\pi l }{2n}\biggr)\sin\phi
  -i\sinh\biggl(\frac{\pi l}{2n}\biggr)\cos\phi\biggr),
\label{eqn:4b}
\end{equation}
which is not zero for any $\phi\in\Sset^1$, as in \eqref{eqn:b2}.

\begin{prop}
\label{prop:4b}
Let $\hat{\Gamma}$ be a domain in $\Cset/(2\pi l/\nu)\Zset$
 containing $\Rset/(2\pi l/\nu)\Zset$ and $\tau=i\alpha K(1/\sqrt{2})$.
The periodically forced Duffing oscillator \eqref{eqn:duf} with $a=0$ 
 is meromorphically nonintegrable in the meaning of Theorem~$\ref{thm:main}$
 near any periodic orbit on the resonant torus
 $\T^\alpha=\{(x^\alpha(\tau),\nu\tau+\theta)\mid\tau\in\hat{\Gamma},\theta\in\Sset^1,\alpha=\alpha^{l/n}\}$
 for $l,n>0$ relatively prime integers.
\end{prop}

\begin{rmk}\
\label{rmk:4b}
\begin{itemize}
\setlength{\leftskip}{-1.8em}
\item[(i)]
As in Remark~{\rm\ref{rmk:4a}(i)},
 if $\beta=0$, then Proposition~$\ref{prop:4b}$ says nothing
 about the nonintegrability of \eqref{eqn:duf}
 since the integral \eqref{eqn:4b} is identically zero.
\item[(ii)]
For any neighborhood $U$ of $\alpha\in(0,\infty)$
 there is not a key set $D\subset U$ for $C^\omega(U)$ such that $M^{l/n}(\phi)$ is not constant
 for $\alpha\in D$ satisfying \eqref{eqn:res4b}.
Hence, Theorem~$\ref{thm:MY}$ is not applicable, as in Remark~{\rm\ref{rmk:4a}(ii)}. 
\end{itemize}
\end{rmk}

\subsection{Case of $a=-1$}

\begin{figure}[t]
\includegraphics[scale=0.26
]{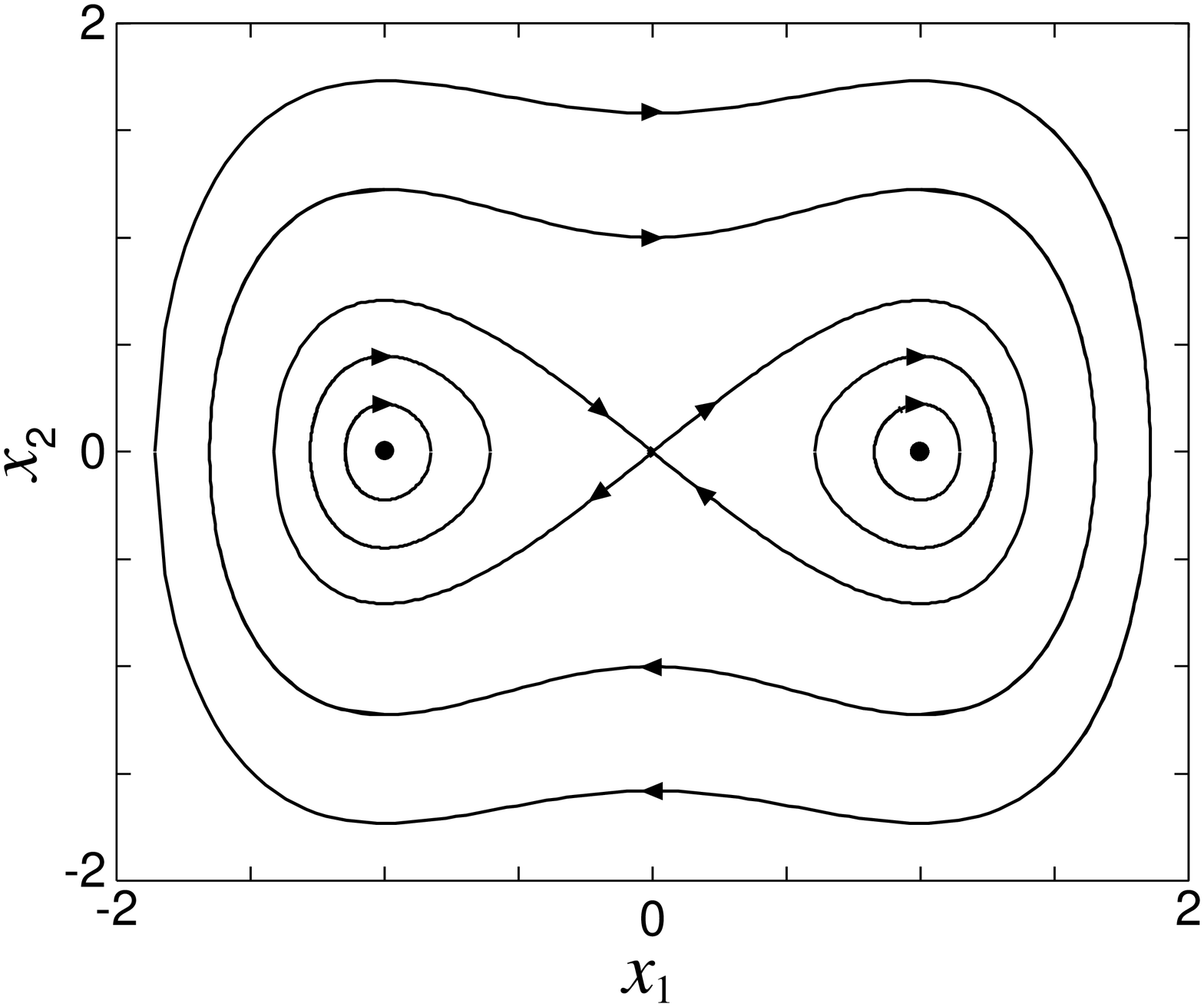}
\caption{Phase portraits of \eqref{eqn:duf} with $\epsilon=0$ and $a=-1$.
\label{fig:4c}}
\end{figure}

We turn to the case of $a=-1$ in \eqref{eqn:duf}.
The phase portraits of \eqref{eqn:duf} with $\epsilon=0$ are shown in Fig.~\ref{fig:4c}.
In particular, there exist a pair of homoclinic orbits
\[
x^\h_\pm(t)=(\pm\sqrt{2}\sech t, \mp\sqrt{2}\sech t\,\tanh t),
\]
a pair of one-parameter families of periodic orbits
\begin{align*}
x^k_{\pm}(t)
 =&\biggl(\pm\frac{\sqrt{2}}{\sqrt{2-k^2}}\dn\left(\frac{t}{\sqrt{2-k^2}}\right),\\
& \quad\mp\frac{\sqrt{2}k^2}{2-k^2}\sn\left(\frac{t}{\sqrt{2-k^2}}\right)
\cn\left(\frac{t}{\sqrt{2-k^2}}\right)\biggr),\quad
k\in(0,1),
\end{align*}
inside each of them, and a one-parameter periodic orbits
\begin{align*}
\tilde{x}^k(t)
 =&\biggl(\frac{\sqrt{2}k}{\sqrt{2k^2-1}}\cn\left(\frac{t}{\sqrt{2k^2-1}}\right),\\
& \quad -\frac{\sqrt{2}k}{2k^2-1}\sn\left(\frac{t}{\sqrt{2k^2-1}}\right) 
 \dn\left(\frac{t}{\sqrt{2k^2-1}}\right)\biggr),\quad
k\in\bigl(1/\sqrt{2},1\bigr),
\end{align*}
outside of them.
The periods of $x^k_{\pm}(t)$ and $\tilde{x}^k(t)$ are given
 by $T^k=2K(k)\sqrt{2-k^2}$ and $\tilde{T}^k=4K(k)\sqrt{2k^2-1}$, respectively
 (see \cite{GH83a,GH83,W90}).
Note that $x^k_{\pm}(t)$ approach $x^\h_{\pm}(t)$ as $k\to 1$.
Assume that the resonance conditions
\begin{equation}
n\hat{T}^k=\frac{\pi l}{\nu},\quad\mbox{i.e.,}\quad
\nu=\frac{\pi l}{nK(k)\sqrt{2-k^2}},
\label{eqn:res4c1}
\end{equation}
and
\begin{equation}
n\tilde{T}^k=\frac{2\pi l}{\nu},\quad\mbox{i.e.,}\quad
\nu=\frac{\pi l}{2nK(k)\sqrt{2k^2-1}},
\label{eqn:res4c2}
\end{equation}
hold at $k=k^{l/n}$ with $l,n>0$ relatively prime integers
 for $x_\pm^k(t)$ and $\tilde{x}^k(t)$, respectively.
We compute the subharmonic Melnikov function \eqref{eqn:mel} as
\[
M_\pm^{l/n}(\tau)=-\delta J_1(k,n)\pm\beta J_2(k,l,n)\sin\tau
\]
and
\[
\tilde{M}^{l/n}(\tau)=-\delta\tilde{J}_1(k,n)+\beta\tilde{J}_2(k,l,n)\sin\tau,
\]
for $x_\pm^k(t)$ and $\tilde{x}^k(t)$, respectively, where
\begin{align*}
&
J_1(k,n)=\frac{4n[(2-k^2)E(k)-2k'^2K(k)]}{3(2-k^2)^{3/2}},\\
&
J_2(k,l,n)=
\begin{cases}
\sqrt{2}\pi \nu\sech\left(\displaystyle\frac{\pi lK(k')}{K(k)}\right) & \mbox{(for $n=1$)};\\
0\quad & \mbox{(for $n\neq 1$)},\
\end{cases}\\
&
\tilde{J}_1(k,n)=\frac{8n[(2k^2-1)E(k)+k'^2K(k)]}{3(2k^2-1)^{3/2}},\\
&
\tilde{J}_2(k,l,n)=
\begin{cases}
2\sqrt{2}\pi \nu\sech\left(\displaystyle\frac{\pi lK(k')}{2K(k)}\right)
 & \mbox{(for $n=1$ and $l$ odd)}; \\
0 & \mbox{(for $n\neq 1$ or $l$ even).}
\end{cases}
\end{align*}
See also \cite{GH83a,GH83,W90} for the computations of the Melnikov functions.

On the other hand, we write the integral \eqref{eqn:thm4a} as
\begin{align}
\hat{\I}(\phi)=&
-\frac{2k^4\delta}{(2-k^2)^2}\int_{\gamma_\phi}
 \sn^2\left(\frac{\tau}{\sqrt{2-k^2}}\right)
 \cn^2\left(\frac{\tau}{\sqrt{2-k^2}}\right)\d\tau\notag\\
&
\mp\frac{\sqrt{2}k^2\beta}{2-k^2}\int_{\gamma_\phi}
 \sn\left(\frac{\tau}{\sqrt{2-k^2}}\right)
  \cn\left(\frac{\tau}{\sqrt{2-k^2}}\right)\cos(\nu\tau+\phi)\d\tau
\label{eqn:c1}
\end{align}
and
\begin{align}
\hat{\I}(\phi)=&
-\frac{2k^2\delta}{(2k^2-1)^2}\int_{\gamma_\phi}
 \sn^2\left(\frac{\tau}{\sqrt{2k^2-1}}\right)
 \dn^2\left(\frac{\tau}{\sqrt{2k^2-1}}\right)\d\tau\notag\\
&
-\frac{\sqrt{2}k\beta}{2k^2-1}\int_{\gamma_\phi}
 \sn\left(\frac{\tau}{\sqrt{2k^2-1}}\right)
  \dn\left(\frac{\tau}{\sqrt{2k^2-1}}\right)\cos(\nu\tau+\phi)\d\tau
\label{eqn:4c1}
\end{align}
for $x_\pm^k(t)$ and $\tilde{x}^k(t)$, respectively.
We take circles centered at $\tau=i\sqrt{2-k^2}K(k')$ and $\tau=i\sqrt{2k^2-1}K(k')$
 with sufficiently small radii as $\gamma_\phi$, and compute \eqref{eqn:c1} and \eqref{eqn:4c1} as
\begin{align}
\hat{\I}(\phi)
=\mp 2\sqrt{2}\pi\nu\beta
 \biggl(\cosh\biggl(\frac{\pi l K(k')}{nK(k)}\biggr)\sin\phi
  &-i\sinh\biggl(\frac{\pi l K(k')}{nK(k)}\biggr)\cos\phi\biggr),
\label{eqn:c2}
\end{align}
and
\begin{align}
\hat{\I}(\phi)
=-2\sqrt{2}\pi\nu\beta
 \biggl(\cosh\biggl(\frac{\pi l K(k')}{2nK(k)}\biggr)\sin\phi
 &-i\sinh\biggl(\frac{\pi l K(k')}{2nK(k)}\biggr)\cos\phi\biggr),
\label{eqn:4c2}
\end{align}
respectively.
See Appendix~C for the derivation of \eqref{eqn:c2}.
The expression \eqref{eqn:4c2} is derived as in \eqref{eqn:b2}.
Note that the integrals \eqref{eqn:c2} and \eqref{eqn:4c2} are not zero for any $\phi\in\Sset^1$.
Applying Theorem~\ref{thm:3a}, we obtain the following.

\begin{prop}
\label{prop:4c}
Let $\hat{\Gamma}$ be a domain in $\Cset/(2\pi l/\nu)\Zset$
 containing $\Rset/(2\pi l/\nu)\Zset$ and $\tau=i\sqrt{2-k^2}K(k')$ $($resp. $\tau=i\sqrt{2k^2-1}K(k'))$.
The periodically forced Duffing oscillator \eqref{eqn:duf} with $a=-1$ 
 is meromorphically nonintegrable in the meaning of Theorem~$\ref{thm:main}$
 near any periodic orbit on the resonant torus
 $\T^k=\{(x^k(\tau),\nu\tau+\theta)\mid\tau\in\hat{\Gamma},\theta\in\Sset^1,k=k^{l/n}\}$
 $($resp. $\T^k=\{(\tilde{x}^k(\tau),\nu\tau+\theta)\mid\tau\in\hat{\Gamma},\theta\in\Sset^1,k=k^{l/n}\})$
 for $l,n>0$ relatively prime integers.
\end{prop}

\begin{rmk}\
\label{rmk:4c}
\begin{itemize}
\setlength{\leftskip}{-1.8em}
\item[(i)]
If $\beta=0$, then Propositions~$\ref{prop:4c}$ says nothing
 about the nonintegrability of \eqref{eqn:duf}
 since the integral \eqref{eqn:b2} is identically zero.
\item[(ii)]
For any neighborhood $U$ of $k=1$
 there is a key set $D\subset U$ for $C^\omega(U)$
 such that $M^{l/n}(\phi)$ $($resp. $\tilde{M}^{l/n}(\phi))$ is not constant
 for $k\in D$ satisfying \eqref{eqn:res4c1} $($resp. \eqref{eqn:res4c2}$)$.
Hence, Theorem~$\ref{thm:MY}$ is applicable to show that
 the periodically forced Duffing oscillator \eqref{eqn:duf} with $a=-1$
 is real-analytic nonintegrable near the surface
 $(\{x^\h(t)\mid t\in\Rset\}\cup\{0\})\times\Sset^1$.
\end{itemize}
\end{rmk}


\section{Additional Examples}

We give two more examples to illustrate Theorem~\ref{thm:main}.

\subsection{Second-order coupled oscillators}
Let $m=\ell$ and consider
\begin{equation}
\dot{I}_j=\epsilon\biggl(-\delta I_j+\Omega_j
 +\beta\sum_{k=1}^\ell\frac{\sin(\theta_k-\theta_j)}{1-\kappa\cos(\theta_k-\theta_j)}\biggr),\quad
\dot{\theta}_j=I_j,\quad
j=1,\ldots,\ell,
\label{eqn:ex2}
\end{equation}
where $\delta,\beta,\kappa,\Omega_j>0$, $j=1,\ldots,\ell$, are constants such that $\kappa<1$.
Equation~\eqref{eqn:ex2} is rewritten in a system of second-order differential equations as
\[
\ddot{\theta}_j+\epsilon\delta\dot{\theta}
 =\epsilon\biggl(\Omega_j
 +\beta\sum_{l=1}^\ell\frac{\sin(\theta_l-\theta_j)}{1-\kappa\cos(\theta_l-\theta_j)}\biggr),\quad
j=1,\ldots,\ell,
\]
which is often referred to as \emph{second-order Kuramoto model} \cite{RPJK16}
 when $\kappa=0$.
When $\delta,\Omega_j=0$, $j=1,\ldots,\ell$,
 the system \eqref{eqn:ex2} is an $\ell$-degree-of-freedom Hamiltonian system
 with the Hamiltonian
\[
H(I,\theta)=\tfrac{1}{2}|I|^2+\frac{\epsilon\beta}{\kappa}
 \sum_{j=2}^\ell\sum_{l=1}^{j-1}\log(1-\kappa\cos(\theta_l-\theta_j)).
\]
Henceforth we only treat a special case of condition (A1) in which
\[
2I_1=I_2=\cdots=I_\ell\neq 0
\]
although infinitely many resonances of multiplicity $\ell-1$ can occur in \eqref{eqn:ex2}.

Let $\omega^\ast=I_1$, so that $T^\ast=2\pi/I_1$.
Assume that
\[
|\theta_j-\theta_k|\neq|\theta_1-\theta_2|\quad\mbox{for $(j,k)\neq(1,2)$},
\] 
and let $\gamma_\theta$ be a closed loop with center at
\[
\tau=\frac{\theta_1-\theta_2}{\omega^\ast}
 +\frac{i}{\omega^\ast}\arccosh\biggl(\frac{1}{\kappa}\biggr)=:\tau^\ast,
\]
and sufficiently small radius.
Using the method of residues,
 we compute the first and second components of \eqref{eqn:A2} as
\[
\I_1(\theta)=\beta\int_{\gamma_\theta}\frac{\sin(\omega^\ast\tau+\theta_2-\theta_1)}
 {1-\kappa\cos(\omega^\ast\tau+\theta_2-\theta_1)}\d\tau
 =2\pi i\kappa\omega^\ast
\]
and
\[
\I_2(\theta)=-\beta\int_{\gamma_\theta}\frac{\sin(\omega^\ast\tau+\theta_2-\theta_1)}
 {1-\kappa\cos(\omega^\ast\tau+\theta_2-\theta_1)}\d\tau
 =-2\pi i\kappa\omega^\ast,
\]
respectively, while its other components are zero. 
Applying Theorem~\ref{thm:main}, we obtain the following.

\begin{prop}
\label{prop:5b}
Let $\Gamma$ be a domain in $\Cset/T^\ast\Zset$
 containing $\Rset/T^\ast\Zset$ and $\tau=\tau^\ast$.
The system \eqref{eqn:ex2} is nonintegrable near {\color{black}any periodic orbit on}
\[
\{(I,\omega^\ast\tau+\theta)\mid
\tau\in\Gamma,I\in\Rset^\ell,\theta\in\Tset^\ell,2I_1=I_2=\cdots=I_\ell\neq 0\}
\]
in the meaning of Theorem~$\ref{thm:main}$.
\end{prop}

\subsection{Pendulum-type oscillator with a constant torque}
We finally set $\ell=m=1$ and consider the two-dimensional system
\begin{equation}
\dot{I}=\epsilon\biggl(\frac{\sin\theta}{1-\kappa\cos\theta}+1\biggr),\quad
\dot{\theta}=I,
\label{eqn:ex1}
\end{equation}
where $\kappa\in(0,1)$ is a constant.
When $\kappa=0$,
 Eq.~\eqref{eqn:ex1} represents an equation of motion
  for the pendulum subjected to a constant torque.
A similar example was treated in \cite{MY21b}.
Assumption~(A1) holds for any $I^\ast=I\neq 0$ as $\omega^\ast=I$ and $T^\ast=2\pi/I$.
Let $\gamma_\theta$ be a closed loop with center at
\[
\tau=-\frac{\theta}{I}+\frac{i}{I}\arccosh\biggl(\frac{1}{\kappa}\biggr)=:\tau^\ast,
\]
and sufficiently small radius, as in Section~5.2.
Noting that $\D\omega(I)=1$ and using the method of residues,
 we compute \eqref{eqn:A2} as
\[
\I(\theta)=\int_{\gamma_\theta}\frac{\sin(I\tau+\theta)}{1-\kappa\cos(I\tau+\theta)}\d\tau
=2\pi i\kappa\omega^\ast.
\]
Applying Theorem~\ref{thm:main}, we obtain the following.

\begin{prop}
\label{prop:5a}
Let $\Gamma$ be a domain in $\Cset/T^\ast\Zset$
 containing $\Rset/T^\ast\Zset$ and $\tau=\tau^\ast$.
The system \eqref{eqn:ex1} is nonintegrable near any periodic orbit
 $\{(I,\omega^\ast\tau+\theta\mid\tau\in\Gamma\}$ for any $I\in\Rset$ and $\theta\in\Sset^1$
 in the meaning of Theorem~$\ref{thm:main}$.
\end{prop}

We easily see that the system \eqref{eqn:ex1} has the first integral
\[
F_1(I,\theta)=\tfrac{1}{2}I^2-\epsilon(\log(1-\kappa\cos\theta)+\theta)
\]
and it is integrable as a system on $\Rset\times\Rset$,
 although $F_1(I,\theta)$ is not even a function on $\Rset\times\Sset^1$.

\section*{Acknowledgements}
The author thanks Shoya Motonaga for helpful and useful discussions.
This work was partially supported by the JSPS KAKENHI Grant Number JP17H02859.


\appendix

\renewcommand{\theequation}{\Alph{section}.\arabic{equation}}
\setcounter{equation}{0}

\section{Previous related results for \eqref{eqn:aasys}}
In this appendix we review some previous related results for the integrability of \eqref{eqn:aasys}.
We begin with the work of Kozlov \cite{K96}.

We first expand $h(I,\theta;0)$ in Fourier series as
\[
h(I,\theta;0)=\sum_{r\in\Zset^m}\hat{h}_r(I)\exp(i r\cdot\theta),
\]
where $\hat{h}_r(I)$, $r\in\Zset^m$, are the Fourier coefficients,
 and assume the following for \eqref{eqn:aasys}:
\begin{enumerate}
\setlength{\leftskip}{-0.8em}
\item[\bf(K1)]
The system~\eqref{eqn:aasys} has $s$ first integrals $F_j(I,\theta)$, $l=1,\ldots,s$,
 which are real-analytic in $(I,\theta,\epsilon)$;
\item[\bf(K2)]
If $r\in\Zset^m$ and $r\cdot\omega(I)=0$ for any $I\in\Rset^\ell$, then $r=0$.
\end{enumerate}
Under assumptions (K1) and (K2)
 we can show that $F_j(I,\theta;0)$, $j=1,\ldots,s$, are independent of $\theta$
  (see Lemma~1 in Section~1 of Chapter IV of \cite{K96}),
 and write $F_{j0}(I;0)=F_j(I,\theta;0)$ and $F_0(I)=(F_{10}(I),\ldots,F_{s0}(I))$.
We refer to $\P_s\subset\Rset^\ell$ as a \emph{Poincar\'e set}
 if for each $I\in\P_s$ there exists linearly independent vectors $r_j\in\Zset^m$, $j=1,\ldots,\ell-s$,
 such that
\begin{enumerate}
\setlength{\leftskip}{-1.8em}
\item[(i)]
$r_j\cdot\omega(I)=0$, $j=1,\ldots ,\ell-s$;
\item[(ii)]
$\hat{h}_{r_j}(I)$, $j=1,\ldots ,\ell-s$, are linearly independent.
\end{enumerate}
Let $U$ be a domain in $\Rset^\ell$.
A set $\Delta\subset U$ is called a \emph{key set} (or \emph{uniqueness set}) for $C^\omega(U)$
 if any analytic function vanishing on $\Delta$ vanishes on $U$.
For example, any dense set in $U$ is a key set for $C^\omega(U)$.
In this situation, we can prove the following theorem
 (see Section~1 of Chapter IV of \cite{K96} for its proof).

\begin{thm}[Kozlov]
\label{thm:Kozlov}
Suppose that assumptions~{\rm(K1)} and {\rm(K2)} hold,
 the Jacobian matrix $\D F_0(I)$ has a maximum rank at a point $I^\ast\in\Rset^\ell$
 and a Poincar\'e set $\P_s\subset U$ is a key set for $C^\omega(U)$,
 where $U$ is a neighborhood of $I^\ast$ in $\Rset^\ell$.
Then the system \eqref{eqn:aasys} has no first integral
 which is real-analytic in $(I,\theta,\epsilon)$
 and functionally independent of $F_j(I,\theta;\epsilon)$, $j=1,\ldots,s$, in $U\times\Tset^m$
 near $\epsilon=0$.
\end{thm}

A version of Theorem~\ref{thm:Kozlov} for the Hamiltonian case $\ell=m$
 was given in \cite{K83} earlier (see also Theorem~7.1 of \cite{AKN06}).
When $s=0$ in (K1),
 Theorem~\ref{thm:Kozlov} means that under the hypotheses
 there exists no first integral which is real-analytic in $(I,\theta,\epsilon)$.
When $s=1$ in (K1), which always occurs if the system \eqref{eqn:aasys} is Hamiltonian,
 it means that under the hypotheses, which hold for $\ell,m=2$
 if besides (K1) and (K2) there exists a key set $\P_1$ for $C^\omega(U)$
 with $\D F_{10}(I)\neq 0$ at a point of $U$
 such that $r\cdot\omega(I)=0$ and $\hat{h}_r(I)\neq 0$ for some $r\in\Zset^2$ on $\P_1$,
 there exists no first integral which is real-analytic in $(I,\theta,\epsilon)$
 and functionally independent of $F_1(I,\theta,\epsilon)$.
In the Hamiltonian case,
 the conclusion implies that the system~\eqref{eqn:aasys} is not Liouville-integrable
 in such a meaning of Theorem~\ref{thm:main}.
However, in the non-Hamiltonian case, this is not generally true:
  it may be Bogoyavlenskij-integrable since it may have $m+\ell-1$ commutative vector fields
  satisfying Definition~\ref{dfn:1a}.
Thus, it is difficult from Theorem~\ref{thm:Kozlov} to say anything
 about Bogoyavlenskij-integrabilty of non-Hamiltonian systems directly. 
 
On the other hand, Motonaga and Yagasaki \cite{MY21b} recently
 discussed nonintegrability of perturbations of general analytically integrable systems
 such that the first integrals and commutative vector fields
 depend analytically on the small parameter, based on the result of \cite{MY21a}.
Let $U$ be a domain in $\Rset^\ell$, as above.
We assume the following:
\begin{enumerate}
\setlength{\leftskip}{0.1em}
\item[\bf(MY1)]
A resonance of multiplicity $m-1$,
\[
\dim_\Qset\langle\omega_1(I),\ldots,\omega_m(I)\rangle=1,
\]
occurs with $\omega(I)\neq 0$ for $I\in D_\R$,
 where $D_\R$ is a key set for $C^\omega(U)$.
\item[\bf(MY2)]
For some $I^\ast\in U$ $\rank\D\omega(I^\ast)=\ell$.
\end{enumerate}
Assumption~(MY1) is similar to assumption~(A1) in Section~1 but more restrictive.
We easily see that if $\rank\D\omega(\bar{I})=m$ for some $\bar{I}\in\Rset^\ell$,
 then assumption~(MY1) as well as (K2) hold for a neighborhood $U$ of $\bar{I}$ in $\Rset^\ell$. 
In (MY1) we take a constant $T_I>0$ for $I\in D_\R$ such that 
\[
\omega_j(I)T_I\in 2\pi\Zset,\quad
j=1,\ldots,m.
\]
Let
\begin{equation}
\bar{\I}_I(\theta)=\int_0^{T_I}h(I,\omega(I)\tau+\theta;0)\d\tau.
\label{eqn:Xi0}
\end{equation}
Their result is stated for \eqref{eqn:aasys} as follows.

\begin{thm}[Motonaga and Yagasaki]
\label{thm:MY}
Suppose that assumptions~{\rm(K2)}, {\rm(MY1)} and {\rm(MY2)} hold.
If there exists a key set $D\subset D_\R$ for $C^\omega(U)$
 such that $\bar{\I}_I(\theta)$ is not constant for $I\in D$,
 then for $|\epsilon|\neq 0$ sufficiently small
 the system \eqref{eqn:aasys} is not real-analytically integrable in the Bogoyavlenskij sense
 in $U\times\Tset^m$ such that
 the first integrals and commutative vector fields
 also depend real-analytically on $\epsilon$ near $\epsilon=0$.
\end{thm}

\begin{rmk}\
\label{rmk:MY}
\begin{enumerate}
\setlength{\leftskip}{-1.8em}
\item[\rm(i)]
If assumption~{\rm(A1)} with $\rank\D\omega(I^\ast)=m$ holds,
 then we can take a neighborhood of the resonant torus $\T^\ast$
 as $U\times\Tset^m$ in Theorem~$\ref{thm:MY}$, like Theorem~$\ref{thm:main}$.
See Section~$2$ of {\rm\cite{MY21a}} for the details.
\item[\rm(ii)]
The integral can also be expressed by the Fourier coefficient $\hat{h}_r(I)$, $r\in\Zset^m$.
See Section~$4$ of {\rm\cite{MY21a}} for the details.
\end{enumerate}
\end{rmk}

Using Theorem~\ref{thm:MY},
 we can discuss Bogoyavlenskij-integrabilty of \eqref{eqn:aasys} even in the non-Hamiltonian case.
However, to determine whether a specific system of the form \eqref{eqn:aasys} is nonintegrable
 in the meaning of Theorem~\ref{thm:MY} or not,
 we need to show that $\bar{\I}_I(\theta)$ is not constant for  infinitely many values of $I$
 since the key set $D$ is an infinite set.
See Section~4 of \cite{MY21b} for more details.

\section{Derivation of \eqref{eqn:b2}}

We use the method of residues and compute the integral \eqref{eqn:b1}.
We begin with the first term in \eqref{eqn:b1}.
Letting $s=1/\sn\zeta$, we have
\begin{equation}
\int\sn^2\zeta\,\dn^2\zeta\,\d\zeta
=-\int\frac{1}{s^4}\sqrt{\frac{k^2-s^2}{1-s^2}}\d s
\label{eqn:b3}
\end{equation}
from the basic properties of the Jacobi elliptic functions,
\begin{equation}
\frac{\d}{\d\zeta}\sn\zeta=\cn\zeta\dn\zeta,\quad
\cn^2\zeta=1-\sn^2\zeta,\quad
\dn^2\zeta=1-k^2\sn^2\zeta.
\label{eqn:b4}
\end{equation}
Obviously, the intgrand in the right hand side of \eqref{eqn:b3} has a pole of order $4$ at $s=0$.
Since $s=1/\sn\zeta=0$ when $\zeta=iK(k')$ and
\[
\frac{\d^3}{\d s^3}\sqrt{\frac{k^2-s^2}{1-s^2}}=0
\]
at $s=0$, we obtain
\begin{equation*}
\int_{\hat{\gamma}_\phi}\sn^2\zeta\,\dn^2\zeta\,\d\zeta
=-\int_{|s|=\rho}\frac{1}{s^4}\sqrt{\frac{k^2-s^2}{1-s^2}}\d s=0
\end{equation*}
by the method of residues,
 where $\hat{\gamma}_\phi=\{\zeta\in\Cset\mid \zeta/\sqrt{1-2k^2}=\gamma_\phi\}$
 and $\rho>0$ is sufficiently small. 
 
We turn to the second term in \eqref{eqn:b1}.
We have
\[
\frac{\d}{\d\zeta}\cn \zeta=-\sn\zeta\,\dn\zeta
=\frac{i}{k(\zeta-iK(k'))^2}+O(1)
\]
near $\zeta=iK(k')$ since
\[
\cn\zeta=-\frac{i}{k(\zeta-iK(k'))}+O(1).
\]
Hence,
\begin{align*}
&
\sn\zeta\,\dn\zeta\,\cos\bigl(\sqrt{1-2k^2}\,\nu\zeta\bigr)\\
&
=-\frac{i\cosh\bigl(\sqrt{1-2k^2}\,\nu K(k')\bigr)}{k(\zeta-iK(k'))^2}
 -\frac{\nu\sqrt{1-2k^2}\sinh\bigl(\sqrt{1-2k^2}\,\nu K(k')\bigr)}{k(\zeta-iK(k'))}+O(1)
\end{align*}
near $\zeta=iK(k')$, so that
\begin{align*}
\int_{\hat{\gamma}_\phi}\sn\zeta\,\dn\zeta\,\cos\bigl(\sqrt{1-2k^2}\,\nu\zeta\bigr)\d\zeta
=&-\frac{2\pi i\nu\sqrt{1-2k^2}}{k}
\sinh\biggl(\frac{\pi l K(k')}{2nK(k)}\biggr),
\end{align*}
where we have used the relation \eqref{eqn:res4a}.
Similarly,
\begin{align*}
\int_{\hat{\gamma}_\phi}\sn\zeta\,\dn\zeta\,\sin\bigl(\sqrt{1-2k^2}\,\nu\zeta\bigr)\d\zeta
 =
 \frac{2\pi i\nu\sqrt{1-2k^2}}{k}\cosh\biggl(\frac{\pi l K(k')}{2nK(k)}\biggr).
\end{align*}
Thus, we obtain \eqref{eqn:b2}.

\section{Derivation of \eqref{eqn:c2}}
We use the method of residues and compute the integral \eqref{eqn:c1}, as in Appendix~B.
We begin with the first term in \eqref{eqn:c1}.
Letting $s=1/\sn\zeta$, we have
\begin{equation}
\int\sn^2\zeta\,\cn^2\zeta\,\d\zeta
=-\int\frac{1}{s^4}\sqrt{\frac{1-s^2}{k^2-s^2}}\d s
\label{eqn:c3}
\end{equation}
by \eqref{eqn:b4}.
Obviously, the intgrand in the right hand side of \eqref{eqn:c3} has a pole of order $4$ at $s=0$.
Since $s=1/\sn\zeta=0$ when $\zeta=iK(k')$ and
\[
\frac{\d^3}{\d s^3}\sqrt{\frac{1-s^2}{k^2-s^2}}=0
\]
at $s=0$, we obtain
\begin{equation*}
\int_{\hat{\gamma}_\phi}\sn^2\zeta\,\cn^2\zeta\,\d\zeta
=-\int_{|s|=\rho}\frac{1}{s^4}\sqrt{\frac{1-s^2}{k^2-s^2}}\d s=0
\end{equation*}
by the method of residues,
 where $\hat{\gamma}_\phi=\{\zeta\in\Cset\mid \zeta/\sqrt{2-k^2}=\gamma_\phi\}$
 and $\rho>0$ is sufficiently small. 

We turn to the second term in \eqref{eqn:c1}.
We have
\[
\frac{\d}{\d\zeta}\dn \zeta=-k^2\sn\zeta\,\cn\zeta
=\frac{i}{(\zeta-iK(k'))^2}+O(1)
\]
near $\zeta=iK(k')$ since
\[
\dn\zeta=-\frac{i}{\zeta-iK(k')}+O(1).
\]
Hence,
\begin{align*}
&
\sn\zeta\,\cn\zeta\,\cos\bigl(\sqrt{2-k^2}\,\nu\zeta\bigr)\\
&
=-\frac{i\cosh\bigl(\sqrt{2-k^2}\,\nu K(k')\bigr)}{k^2(\zeta-iK(k'))^2}
 -\frac{\nu\sqrt{2-k^2}\sinh\bigl(\sqrt{2-k^2}\,\nu K(k')\bigr)}{k^2(\zeta-iK(k'))}+O(1)
\end{align*}
near $\zeta=iK(k')$, so that
\begin{align*}
\int_{\hat{\gamma}_\phi}\sn\zeta\,\cn\zeta\,\cos\bigl(\sqrt{2-k^2}\,\nu\zeta\bigr)\d\zeta
=&-\frac{2\pi i\nu\sqrt{2-k^2}}{k^2}
 \sinh\biggl(\frac{\pi l K(k')}{nK(k)}\biggr),
\end{align*}
where we have used the relation \eqref{eqn:res4c1}.
Similarly,
\begin{align*}
\int_{\hat{\gamma}_\phi}\sn\zeta\,\cn\zeta\,\sin\bigl(\sqrt{2-k^2}\,\nu\zeta\bigr)\d\zeta
 =
 \frac{2\pi i\nu\sqrt{2-k^2}}{k^2}\cosh\biggl(\frac{\pi l K(k')}{nK(k)}\biggr).
\end{align*}
Thus, we obtain \eqref{eqn:c2}.



\begin{thebibliography}{99}

\bibitem{A89}
V.I.~Arnold,
\textit{Mathematical Methods of Classical Mechanics}, 2nd ed.,
Springer, New York, 1989.
\bibitem{AKN06}
V. I.~Arnold, V.V.~Kozlov and A.I.~Neishtadt,
\textit{Dynamical Systems III: Mathematical Aspects of Classical and Celestial Mechanics}, 3rd ed.,
Springer, Berlin, 2006.
\bibitem{AZ10}
M.~Ayoul and N.T.~Zung,
Galoisian obstructions to non-Hamiltonian integrability,
\textit{C. R. Math. Acad. Sci. Paris}, \textbf{348} (2010), 1323--1326.
\bibitem{B96}
J.~Barrow-Green,
\textit{Poincar\'e and the Three-Body Problem},
American Mathematical Society, Providence, RI, 1996.
\bibitem{B98}
O.I.~Bogoyavlenskij, 
Extended integrability and bi-hamiltonian systems, 
\textit{Comm. Math. Phys.}, \textbf{196} (1998), 19--51.
\bibitem{BF54}
P.F.~Byrd and M.D.~Friedman,
\textit{Handbook of Elliptic Integrals for Engineers and Physicists},
Springer, Berlin, 1954.


\bibitem{D18}
G.~Duffing,
\textit{Erzwungene Schwingungen bei Ver\"anderlicher Eigenfrequenz und Ihre Technische Bedeutung},
Sammlung Vieweg, Braunschweig, 1918.



\bibitem{GH83a}
B.D.~Greenspan and P.J.~Holmes,
Homoclinic orbits, subharmonics and global bifurcations in forced oscillations,
in \textit{Nonlinear Dynamics and Turbulence}, G.I.~Barenblatt, G.~Iooss and D.D.~Joseph (eds.) 
Pitman, Boston, MA, 1983.

\bibitem{GH83}
J.~Guckenheimer and P.J.~Holmes,
\textit{Nonlinear Oscillations, Dynamical Systems, and Bifurcations of Vector Fields},
Springer, New York, 1983.

\bibitem{H79}
P.J.~Holmes,
A nonlinear oscillator with a strange attractor,
\textit{ Philos. Trans. Roy. Soc. London Ser. A},
\textbf{292} (1979), 419--448.

\bibitem{K83}
V.V.~Kozlov,
Integrability and non-integarbility in Hamiltonian mechanics,
\textit{Russian Math. Surveys}, \textbf{38} (1983), 1--76.
\bibitem{K96}
V.V.~Kozlov,
\textit{Symmetries, Topology and Resonances in Hamiltonian Mechanics},
Springer, Berlin, 1996.
\bibitem{M63}
V.K.~Melnikov,
On the stability of the center for time periodic perturbations,
\textit{Trans. Moscow Math. Soc.}, \textbf{12} (1963), 1--56.
\bibitem{M99}
J.J.~Morales-Ruiz,
\textit{Differential Galois Theory and Non-Integrability of Hamiltonian Systems},
Birkh\"auser, Basel, 1999.
\bibitem{M02}
J.J.~Morales-Ruiz,
A note on a connection between the Poincar\'e-Arnold-Melnikov integral
 and the Picard-Vessiot theory, 
in \textit{Differential Galois theory}, T.~Crespo and Z.~Hajto (eds.),
Banach Center Publ. 58, 
Polish Acad. Sci. Inst. Math., 2002, pp.~165--175.
\bibitem{MR01}
J.J.~Morales-Ruiz and J.P.~Ramis,
Galoisian obstructions to integrability of Hamiltonian systems,
\textit{Methods, Appl. Anal.}, \textbf{8} (2001), 33--96.
\bibitem{MRS07}
J.J.~Morales-Ruiz, J.-P.~Ramis and C.~Simo,
Integrability of Hamiltonian systems and differential Galois groups of higher variational equations, \textit{Ann. Sci. \'Ecole Norm. Suppl.}, \textbf{40} (2007), 845--884.

\bibitem{M73}
J.~Moser,
\textit{Stable and Random Motions in Dynamical Systems},
Princeton University Press, Princeton, 1973.

\bibitem{MY18}
S.~Motonaga and K.~Yagasaki,
Nonintegrability of parametrically forced nonlinear oscillators,
\textit{Regul. Chaotic Dyn.}, \textbf{23} (2018), 291--303.
\bibitem{MY21a}
S.~Motonaga and K.~Yagasaki,
Persistence of periodic and homoclinic orbits, first integrals and commutative vector fields
 in dynamical systems,
\textit{Nonlinearity}, \textbf{34} (2021), 7574--7608.
\bibitem{MY21b}
S.~Motonaga and K.~Yagasaki,
Obstructions to integrability of nearly integrable dynamical systems near regular level sets, 
submitted for publication.

\bibitem{P90}
H.~Poincar\'e,
Sur le probl\'eme des trois corps et les \'equations de la dynamique,
\textit{Acta Math.}, \textbf{13} (1890), 1--270;
English translation:
\textit{The Three-Body Problem and the Equations of Dynamics},
Translated by D.~Popp, Springer, Cham, Switzerland, 2017.
\bibitem{P92}
H.~Poincar\'e,
\textit{New Methods of Celestial Mechanics}, Vol.~1,
AIP Press,  New York, 1992 (original 1892).


\bibitem{RPJK16}
F.A.~Rodrigues, T.K.D.~Peron, P.~Ji and J.~Kurths,
The Kuramoto model in complex networks,
\textit{Phys. Rep.}, \textbf{610} (2016), 1--98.

%

\bibitem{U78}
Y.~Ueda,
Random phenomena resulting from nonlinearity in the system described by Duffing's equation,
\textit{Internat. J. Non-Linear Mech.}, \textbf{20} (1985), 481--491 (original 1978).

\bibitem{WW27}
E.T.~Whittaker and G.N.~Watson,
\textit{A Course in Modern Analysis}, 4th ed..
Cambridge University Press, Cambridge, 1927.


\bibitem{W90}
S.~Wiggins,
\textit{Introduction to Applied Nonlinear Dynamical Systems and Chaos}, 
Springer, New York, 
1990.
\bibitem{Y94}
K.~Yagasaki,
Homoclinic motions and chaos in the quasiperiodically forced van der Pol-Duffing oscillator
 with single well potential,
 \textit{Proc. R. Soc. Lond. A}, 445 (1994), 597--617.
\bibitem{Y96}
K.~Yagasaki,
The Melnikov theory for subharmonics and their bifurcations in forced oscillations,
\textit{SIAM J. Appl. Math.}, \textbf{56} (1996),1720--1765.
\bibitem{Y02}
K.~Yagasaki,
Melnikov's method and codimension-two bifurcations in forced oscillations,
\textit{J. Differential Equations}, \textbf{185} (2002), 
 1--24.
\bibitem{Y03}
K.~Yagasaki,
Degenerate resonances in forced oscillators,
\textit{Discrete Contin. Dyn. Syst. B}, \textbf{3} (2003), 
\bibitem{Y21a}
K.~Yagasaki,
Nonintegrability of the restricted three-body problem,
submitted for publication.
\bibitem{Y21b}
K.~Yagasaki,
New proof of Poincar\'e's result on the restricted three-body problem,
submitted for publication.
\bibitem{Z82}
S.L.~Ziglin,
Self-intersection of the complex separatrices and the non-existing of the integrals
 in the Hamiltonian systems with one-and-half degrees of freedom,
\textit{J. Appl. Math. Mech.}, \textbf{45} (1982), 411--413.
\bibitem{Z18}
N.T.~Zung,
A conceptual approach to the problem of action-angle variables,
\textit{Arch. Ration. Mech. Anal.}, \textbf{229} (2018), 
789--833.
\end{thebibliography}
\end{document}